\documentclass[11pt]{article}
\usepackage{amsfonts}
\usepackage{a4wide,latexsym,amssymb,amsfonts,amsmath,multirow}
\usepackage{txfonts}
\newtheorem{theorem}{Theorem}[section]

\newtheorem{lemma}[theorem]{Lemma}
\newtheorem{corollary}[theorem]{Corollary}
\newtheorem{example}{Example}[section]
\newtheorem{construction}[theorem]{Construction}

\def\whitebox{{\hbox{\hskip 1pt
\vrule height 6pt depth 1.5pt \lower 1.5pt\vbox to 7.5pt{\hrule
width 3.2pt\vfill\hrule width 3.2pt} \vrule height 6pt depth 1.5pt
\hskip 1pt } }}
\def\qed{\ifhmode\allowbreak\else\nobreak\fi\hfill\quad\nobreak
\whitebox\medbreak}

\newcommand{\pf}{\noindent{\bf Proof.}\ }
\newcommand{\ignore}[1]{}

\begin{document}
\title{\bf Locating arrays with mixed alphabet sizes
\thanks{Correspondence to: Ce Shi (shice060@lixin.edu.cn).
The second author's work was supported by NSFC grant 11301342 and Natural Science Foundation of Shanghai No. 17ZR1419900}}
\author{ Ce Shi$^1$,  Hao Jin$^2$ and Tatsuhiro Tsuchiya $^2$\\
\small $^1$  School of  Statistics and Mathematics \\
\small Shanghai Lixin University of Accounting and Finance, \
\small Shanghai 201209, China \\
\small $^2$ Graduate School of Information Science and Technology\\
\small Osaka University,\
\small Suita 565-0871, Japan \\}

\maketitle
\begin{abstract}
Locating arrays (LAs) can be used to detect and identify interaction faults among factors in a component-based
system. The optimality and constructions of LAs with a single fault have been
investigated extensively under the assumption that all the factors have the same values.
However, in real life, different factors in a system have different numbers of possible values. Thus, it is necessary for LAs to
satisfy such requirements. We herein establish a general lower bound on
the size of mixed-level $(\bar{1},t)$-locating arrays. Some methods for constructing LAs including
direct and recursive constructions are provided. In particular, constructions that
produce optimal LAs satisfying the lower bound are described. Additionally,
some series of optimal LAs satisfying the lower bound are presented.


\end{abstract}

\section{Introduction}
Testing is important in detecting failures triggered by interactions among
factors. As reported in \cite{CMMSSY2006}, owing to the complexity
of information systems, interactions among components are complex and numerous. Ideally, one
would test all possible interactions (exhaustive testing); however, this is often infeasible
owing to the time and cost of tests, even for a moderately small system. Therefore, test suites that provide coverage of the most prevalent interactions should be developed.
Testing strategies that use such test suites are usually called combinatorial testing
or combinatorial interaction testing (CIT). CIT has shown its effectiveness in detecting faults,
particularly in component-based systems or configurable systems \cite{KKL2013,NL2011}.

The primary combinatorial object used to generate a test suite for CIT is covering arrays (CAs).
CAs are applied in the testing of networks, software, and hardware, as well as
construction and related applications \cite{Colbourn2004,KC2019,Sloane1993}. In a CA, the factors have the same number of values. However, in real life,
different factors have different numbers of possible values. Thus, mixed-level CAs
or mixed covering arrays (MCAs) are a natural extension of covering array research, which improves their suitability
for applications \cite{CDPP1996,CDFP1997,CMMSSY2006,CSWY2011,MSSW2003,Sherwood2008}. A CA or MCA
as a test suite can be used to detect the presence of failure-triggered interactions. However, they do not
guarantee that faulty interactions can be identified. Consequently, tests to reveal the location of interaction
faults are of interest. To address this problem, Colbourn and McClary formalized the problem of non-adaptive
location of interaction faults and proposed the notion of locating arrays (LAs) \cite{CM2008}.

LAs are a variant of CAs with the ability to determine faulty interactions from the outcomes of the tests.
An LA with parameters $d$ and $t$ is denoted by $(d,t)$-LA, where $d$ and $t$ represent the numbers of
faulty interactions and of components or factors in a faulty interaction, respectively.
$t$ is often called {\em strength}. When the number of faulty interactions is at most, instead of exactly $d$, we use the
notation $(\bar{d},t)$-LA to denote it. Generally, testing with a $(d,t)$-LA can not only
detect the presence of faulty interactions, but can also identify $d$ faulty interactions. Similarly,
using a $(\bar{d},t)$-LA as a test suite allows one to identify all faulty interactions if the number is at most $d$.

LAs have been utilized in measurement and testing \cite{ACS2015,CMCPS2016,CS2016}.
Mart\'{\i}nez et al. \cite{MMPS} developed adaptive analogues and established feasibility
conditions for an LA to exist. Only the minimum number of tests in $(1,1)$-LA and $(\bar{1},1)$-LA
is known precisely \cite{CFH2017}. The minimum number of rows in an LA is determined when the number
of factors is small \cite{STY2012,TCY2012}.
When $(d,t)=(1,2)$, three recursive constructions are provided, as in \cite{CF2016}. Beyond these few direct
and recursive constructions, computation methods are applied to construct $(1,2)$-LAs using a Constraint Satisfaction Problem (CSP) solver
and a Satisfiability (SAT) solver \cite{KKNT2017,KKNT2019,NKNT2014}. Lanus {\em et al. }\cite{LCM2019} described a
randomized computational search algorithm called partitioned search with column resampling to construct $(1,t)$-LAs.
Furthermore, column resampling can be applied to construct $(\bar{1},t)$-LA with $\delta\leq 4$ \cite{SSCS2018}.
The first and third authors extended the notion of LAs to expand the applicability to practical testing problems. Specifically, they
proposed constrained locating arrays (CLAs), that can be used to detect and locate failure-triggering interactions in the
presence of constraints. Computational constructions for this variant of LAs can be found in \cite{JKCT2018,JT2018,JT20183}.

Although a few constructions exist for $(1,t)$-LAs and $(\bar{1},t)$-LAs, these methods do not treat
cases where different factors have difference values. For real-world applications, it is desirable for LAs to
satisfy such requirements. Herein, we will focus on mixed-level $(\bar{1},t)$-LAs, which is equivalent to mixed-level $(1,t)$-LAs,
and an MCA by Lemma \ref{(1,t)-LA and MCA}.

The remainder of the paper is organized as follows. The next section provides the definitions of basic concepts,
such as MCAs and LAs. A general lower bound on the size of mixed-level $(\bar{1},t)$-LAs
will be established in Section 3, which will be regarded as benchmarks for
the construction of optimal LAs with specific parameters. Some methods for constructing LAs including
direct and recursive constructions are provided in Section 4. In particular, some constructions that produce optimal LAs satisfying the lower bound will
be described in this section. The final section contains some concluding remarks.

\section{Definitions and Notations}

The notation $I_n$ represents the set $\{1,2,\cdots, n\}$, while the notations $N,k$ and $t$ represent positive integers with $t<k$. We herein model CIT as follows. Suppose that $k$ factors
denoted by $F_1,F_2,\cdots, F_k$ exist. The $i$th factor has a set of $v_i$ possible values (levels) from a set $V_i$,
where $i\in I_k$. A test is a $k$-tuple $(a_1,a_2,\cdots, a_k)$, where $a_i\in V_i$ for $1\leq i\leq k$.
A test, when executed, has the following outcome: {\em pass} or {\em fail}. A test suite is a collection of tests,
and the outcomes are the corresponding set of pass/fail results. A fault is evidenced by a failure outcome for a test.
Tests are considered to be executed in parallel; therefore, testing is non-adaptive or predetermined.

Let $A=(a_{ij}) (i \in I_N, j\in I_k)$ be an $N\times k$ array with entries in the $j$th column
from a set $V_j$ of $v_j$ symbols. A {\em $t$-way interaction} is a possible $t$-tuple of values for any $t$-set of columns,
 denoted by $T=\{(i, \sigma_i): \sigma_i\in V_i, i\in I\subseteq I_k, |I|=t\}$.
We denote $\rho (A,T)=\{r: a_{ri}=\sigma_i, i\in I\subseteq I_k, |I|=t\}$ for the set of rows of $A$, in which the interaction
is included. For an arbitrary set $\mathcal {T}$ of $t$-way interactions, we define $\rho (A,\mathcal {T})=\cup_{T\in {\cal T}}\rho (A,T)$.
We use the notation $\mathcal {I}_t$ to denote the set of all $t$-way interactions of $A$.

The array $A$ is termed MCAs, denoted by
MCA$_\lambda(N; t, k, (v_1,v_2,\cdots, v_k))$ if $|\rho (A,T)|\geq \lambda$ for
all $t$-way interactions $T$ of $A$. In other words, $A$ is an MCA if each $N \times t$ sub-array
includes all the $t$-tuples $\lambda$ times at the least. Here, the number of rows $N$ is called the array size. The number $\lambda$ is termed as the array index. The number of columns $k$ is called
the number of factors (or variables), number of components, or degree. The word ``strength'' is generally
accepted for referring to the parameter $t$. When $\lambda=1$, the notation MCA$(N; t, k, (v_1, v_2,\cdots, v_k))$
is used.

When $v_1=v_2=\cdots =v_k=v$, an MCA$_\lambda(N; t, k, (v_1,v_2,\cdots, v_k))$ is merely a CA$_\lambda(N; t,k,v)$.
When $\lambda=1$ in a CA, we omit the subscript. Without loss of generality, we often assume that the symbol set sizes are in a non-decreasing order, i.e., $v_1\leq v_2\leq \cdots \leq v_k$. Hereinafter, these assumptions will continue to
be used. When $v_i=1$, the presence of the $i$th factor does not affect the properties of
the mixed covering arrays; thus, it is often assumed that $v_i\geq 2$ for $1\leq i\leq k$.

Following \cite{CM2008}, if, for any
${\cal T}_1, {\cal T}_2 \subseteq \mathcal {I}_t$ with $|{\cal T}_1|
= |{\cal T}_2| = d$, we have
$$\rho(A, {\cal T}_1) = \rho(A, {\cal T}_2)
\Leftrightarrow {\cal T}_1 = {\cal T}_2,$$

\noindent then the array $A$ is regarded as a $(d,t)$-LA
and denoted by $(d,t)$-LA$(N;k, (v_1,v_2,\cdots, v_k))$. Similarly, the
definition is extended to permit sets of $d$ interactions at the most by
writing $\bar{d}$ in place of $d$ and permitting instead
$|{\cal T}_1| \le d$ and $|{\cal T}_2| \le d$. In this case, we use the notation
$(\bar{d},t)$-LA$(N;k, (v_1,v_2,\cdots, v_k))$. Clearly, the condition $\rho(A, {\cal T}_1) = \rho(A, {\cal T}_2)
\Leftrightarrow {\cal T}_1 = {\cal T}_2$ is satisfied if ${\cal T}_1 \not= {\cal T}_2\Rightarrow \rho(A, {\cal T}_1) \not= \rho(A, {\cal
T}_2).$ In the following, we will fully apply this fact.

We herein focus on $(\bar{1},t)$-LA$(N;k, (v_1,v_2,\cdots, v_k))$ in this paper. One of the main problems
regarding $(\bar{1},t)$-LA$(N; k, (v_1,v_2,\cdots, v_k))$ is the construction of such LAs having the minimum $N$ when
its other parameters have been fixed. However, this is a difficult and challenging problem. The larger the
strength $t$, the more difficult it is to construct a minimum LA. We use the notations
$(\bar{1},t)$-LAN$(k, (v_1,v_2,\cdots, v_k))$ to represent the minimum number $N$, for which a $(\bar{1},t)$-LA$(N;k, (v_1,v_2,\cdots, v_k))$
exists. A $(\bar{1},t)$-LA$(N;k, (v_1,v_2,\cdots, v_k))$ is called {\em optimal} if $N=(\bar{1},t)\mbox{-LAN}(k, (v_1,v_2,\cdots, v_k))$.

\begin{lemma}\cite{KKNT2019}\label{(1,t)-LA and MCA}
Suppose that $A$ is an $N\times k$ array. $A$ is a $(\bar{1},t)$-LA$(N;k, (v_1,v_2,\cdots, v_k))$ if and only if it is
a $(1,t)$-LA$(N;k, (v_1,v_2,\cdots, v_k))$ and an MCA.
\end{lemma}

Lemma \ref{(1,t)-LA and MCA} shows that $A$ is a $(\bar{1},t)$-LA if $A$ is an MCA and $\rho(A,T_1)\not= \rho(A,T_2)$ whenever $T_1$ and $T_2$ are distinct $t$-way interactions. We will use this simple fact hereinafter.

\section{A lower bound on the size of $(\bar{1},t)$-LA$(N;k, (v_1,v_2,\cdots, v_k))$}

A benchmark to measure the optimality for $(\bar{1},t)$-LA$(N;k, (v_1,v_2,\cdots, v_k))$ is described in this section.
It follows from Lemma \ref{(1,t)-LA and MCA} that $A$ is a $(\bar{1},t)$-LA only if $A$ is an MCA,
which implies that $|\rho(A,T)|\geq 1$ for any $t$-way interaction $T$ of $A$.
Consequently, $(\bar{1},t)$-LAN$(k, (v_1,v_2,\cdots, v_k))\geq \prod _{i=k-t+1}^k v_i$, where $2\leq v_1\leq v_2\leq \cdots \leq v_k$.
Specifically, we have the following results.

\begin{lemma}\label{case 1}
Let $2\leq v_1\leq v_2\leq \cdots \leq v_{k-t}, 2v_{k-t}\leq v_{k-t+1}\leq \cdots \leq v_k$.
Then, $(\bar{1},t)$-LAN$(k, (v_1,v_2,\cdots,v_k))\\ \geq \prod _{i=k-t+1}^k v_i$.
\end{lemma}

It is remarkable that the lower bound on the size of $(\bar{1},t)$-LA$(N;k, (v_1,v_2,\cdots, v_k))$
in Lemma \ref{case 1} can be achieved. We will present some infinite classes of optimal $(\bar{1},t)$-LA$(N;k, (v_1,v_2,\cdots, v_k))$
satisfying the lower bound in the next section. When $v_i=v_{i+1}=\cdots= v_{k-t}= v_{k-t+1}$, where $i\in \{1,2,\cdots, k-t\}$,
we can obtain a lower bound on the size of $(\bar{1},t)$-LA  by the similar argument as the proof of
Theorem 3.1 in \cite{TCY2012}. We state it as follows.

\begin{lemma}\label{case 2}
Let $2\leq v_1\leq v_2\leq \cdots  \leq v_k$. If $v_i=v_{i+1}=\cdots= v_{k-t}= v_{k-t+1}$, where $i\in \{1,2,\cdots, k-t\}$,
then $(\bar{1},t)$-LAN$(k, (v_1,v_2,\cdots, v_k))\geq  \left\lceil\frac{2\sum_{i\leq j_1<\cdots<j_t \le k}
\prod_{s=1}^{t} v_{j_s} }{1+\binom {k-i+1} {t}}\right\rceil$.
\end{lemma}

\pf Let $A$ be a $(\bar{1},t)$-LA$(N;k, (v_1,v_2,\cdots, v_k))$. We can obtain an $N\times (k-i+1)$ array $A'$
by selecting the last $(k-i+1)$ columns of $A$ (if $i=1$, then $A'$ is merely $A$). 
In the array $A'$, for any $i\leq j_1<\cdots<j_t \le k$, we write $n^{\ell}_{j_1\dots j_t} = |S^{\ell}_{j_1\dots j_t}|$, where
$S^{\ell}_{j_1\dots j_t} = \left\{((j_1,x_1), \ldots, (j_t,x_t))\big||\rho(A',((j_1,x_1),  \cdots, (j_t,x_t)))|={\ell}\right\}, \\{\ell}=1,2,3,\dots $.

As stated above, $|\rho(A,T)|\geq 1$ for any $t$-way interaction $T$ of $A'$. Consequently,
$\sum_{{\ell}\ge 1}n^{\ell}_{j_1\dots j_t} = \prod_{s=1}^{t} v_{j_s}$ and $\sum_{{\ell}\ge 1}({\ell} \times n^{\ell}_{j_1\dots j_t}) = N$ hold.
It is deduced that $n^1_{j_1\dots j_t} \ge  2\prod_{s=1}^{t} v_{j_s}-N$.
By Lemma \ref{(1,t)-LA and MCA} and the proof of Lemma \ref{trun}, $A'$ is a $(1,t)$-LA. Thus, in any two of $\binom{k-i+1}{t}$ sets, $\rho(A',S^1_{j_1\dots j_t})'$s
with $i\le j_1<\cdots<j_t \le k$ share no common elements. Hence, $\sum_{i\leq j_1<\cdots<j_t \le k}n^1_{j_1\dots j_t}\leq N$,
which implies that $\sum_{i\leq j_1<\cdots<j_t \le k}(2\prod_{s=1}^{t} v_{j_s} - N )\leq \sum_{i\leq j_1<\cdots<j_t \le k}n^1_{j_1\dots j_t} \leq N$, i.e.,
$N\geq \left\lceil\frac{2\sum_{i\leq j_1<\cdots<j_t \le k} \prod_{s=1}^{t} v_{j_s} }{1+\binom {k-i+1}{t}}\right\rceil$.
Hence, $(\bar{1},t)\mbox{-LAN}(k, (v_1,v_2,\cdots, v_k))\geq \left\lceil\frac{2\sum_{i\leq j_1<\cdots<j_t \le k} \prod_{s=1}^{t} v_{j_s} }{1+\binom {k-i+1}{t}}\right\rceil $. \qed

Based on $i=1$ and $v_{k-t+1}=\cdots=v_k=v$ in Lemma \ref{case 2}, the following corollary can be easily obtained.
It serves as a benchmark for a $(1,t)$-LA$(N;k,v)$, which was first presented in \cite{TCY2012}.

\begin{corollary}
Let $v, t$, and $k$ be integers with $t<k$. Then, $(1,t)-\mbox{LAN}\ (t, k, v) \ge \left\lceil\frac{2\left(^k_t\right) v^t}{1+\left(^k_t\right)}\right\rceil$.
\end{corollary}

In a $(\bar{1},t)$-LA$(N;k, (v_1,v_2,\cdots, v_k))$, we often assume that $2\leq v_1\leq v_2\leq \cdots \leq v_{k-t}\leq v_{k-t+1} \leq \cdots \leq v_k$.
Lemma \ref{case 1} and Lemma \ref{case 2} consider the cases $v_{k-t}=v_{k-t+1}$ and $2v_{k-t}\leq v_{k-t+1}$, respectively. The left case is $v_{k-t}<v_{k-t+1}<2v_{k-t}$,
which is considered in the following lemma.

\begin{lemma}\label{case 3}
Let $2\leq v_1\leq v_2\leq \cdots  \leq v_k$. If $v_{k-t}< v_{k-t+1}<2v_{k-t}$,
then $(\bar{1},t)$-LAN$(k, (v_1,v_2,\cdots, v_k))\geq  m $, where
\begin{eqnarray*}
m = \left\{
\begin{array}{ll}
  \mbox{max}\{\left\lceil\frac{2\sum_{k-t\leq j_1<\cdots<j_t \le k}
\prod_{s=1}^{t} v_{j_s} }{t+2}\right\rceil, \prod_{i=k-t+1}^k v_i+\prod_{i=k-t+2}^k v_i \}, & \mbox{if }\  t\geq 2;\\
\left\lceil\frac{2 v_{k-1}+2v_k }{3}\right\rceil, & \mbox{if}\  t=1.
\end{array}
\right.
\end{eqnarray*}
\end{lemma}
\pf From the above argument, it is known that $(\bar{1},t)$-LAN$(k, (v_1,v_2,\cdots, v_k))\geq M=\prod_{i=k-t+1}^k v_i$.
Suppose that $A$ is a $(\bar{1},t)$-LA$(N;k, (v_1,v_2,\cdots, v_k))$, where $N=M+L$ and $L\geq 0$.
Select the last $(t+1)$ columns of $A$ to form an $N\times (t+1)$ array $A'$. By Lemma \ref{trun}, $A'$ is a $(\bar{1},t)$-LA$(N;t+1,(v_t,v_{t+1},\cdots,v_k))$.
Similar to the proof of Lemma \ref{case 2}, we can prove that $N\geq \left\lceil\frac{2\sum_{k-t\leq j_1<\cdots<j_t \le k} \prod_{s=1}^{t} v_{j_s} }{t+2}\right\rceil$.
When $t=1$, we can obtain $m=\left\lceil\frac{2 v_{k-1}+2v_k }{3}\right\rceil$. For $t\geq 2$, we will prove that $N\geq M+\prod_{i=k-t+2}^k v_i$, i.e., $L\geq \prod_{i=k-t+2}^k v_i$. Without loss of generality, suppose that $A'$ contains two parts, the first part is an $M\times (t+1)$ array $B$ containing an $M\times t$ sub-array comprising all $t$-tuples
over $V_{k-t+1}\times V_{k-t+2}\times \cdots \times V_k$; the left part is an $L\times (t+1)$ array $C$. (If $L=0$, then $B=A'$).

If $L<\prod_{i=k-t+2}^k v_i$, then at least one $(t-1)$-way
interaction $T=\{(i,a_i): i\in I_k\setminus I_{k-t+1},a_i\in V_i\}$ exists such that it is not included by any
row of $C$ (If $B=A'$, then all the $(t-1)$-way interactions satisfy the condition. We can choose an arbitrary one). Hence, we have
$|\rho(A', T_1)|=1$ for any $t$-way interaction $T_1\in \mathcal{T}_1=\{T\cup(k-t+1,i):i\in V_{k-t+1}\}$.
Since $A$ is a $(\bar{1},t)$-LA$(N;k, (v_1,v_2,\cdots, v_k))$,
$|\rho(A',T_2)|\geq 1$ for any $t$-way interaction $T_2\in \mathcal{T}_2=\{T\cup(k-t,i):i\in V_{k-t}\}$.
It is clear that $\rho(A',\mathcal{T}_1)= \rho(B,T)=\rho(A',T)=\rho(A',\mathcal{T}_2)$ with $|\rho(A',\mathcal{T}_1)|=v_{k-t+1}$.

Because $|\mathcal{T}_2|=v_{k-t}<|\mathcal{T}_1|=v_{k-t+1}<2|\mathcal{T}_2|$, at least
one $t$-way interaction $T'\in \mathcal{T}_2$ exists such that $|\rho(A',T')|=1$. Otherwise, $|\rho(A',T')|\geq 2$ for
any $t$-way interaction $T'\in \mathcal{T}_2$, which implies that $|\rho(A',\mathcal{T}_2)|\geq 2|\mathcal{T}_2|=2v_{k-t}$,
but $|\rho(A',\mathcal{T}_2)|=|\rho(A',\mathcal{T}_1)|=v_{k-t+1}<2v_{k-t}$. It follows that $\rho(A',T')=\rho(A',T_1')$,
where $T_1'$ is a certain $t$-way interaction of $\mathcal{T}_1$. It is obvious that $T'\not =T_1'$.
Consequently, $A'$ is not a $(1,t)$-LA. Thus, $L\geq \prod_{i=k-t+2}^k v_i$. Consequently,  $m=\mbox{max}\{\left\lceil\frac{2\sum_{k-t\leq j_1<\cdots<j_t \le k}
\prod_{s=1}^{t} v_{j_s} }{t+2}\right\rceil, \prod_{i=k-t+1}^k v_i+\prod_{i=k-t+2}^k v_i \}$ if $t\geq 2$. \qed

Combining Lemmas \ref{case 1}, \ref{case 2}, and \ref{case 3}, a lower bound on the size
of $(\bar{1},t)$-LA$(N;k, (v_1,v_2,\cdots, v_k))$ can be obtained, which serves as a benchmark
to measure the optimality.

\begin{theorem}\label{L-bound}
Let $2\leq v_1\leq v_2\leq \cdots  \leq v_k$. Then, $(\bar{1},t)$-LAN$(k, (v_1,v_2,\cdots, v_k))\geq$

\begin{enumerate}
\item $\prod _{i=k-t+1}^k v_i$, if $2v_{k-t}\leq v_{k-t+1}$;
\item    $\left\lceil\frac{2\sum_{i\leq j_1<\cdots<j_t \le k}
\prod_{s=1}^{t} v_{j_s} }{1+\binom {k-i+1} {t}}\right\rceil$,                      if $v_i=v_{i+1}=\cdots= v_{k-t}= v_{k-t+1}$, where $i\in \{1,2,\cdots, k-t\}$;

\item $\mbox{max}\{\left\lceil\frac{2\sum_{k-t\leq j_1<\cdots<j_t \le k}
\prod_{s=1}^{t} v_{j_s} }{t+2}\right\rceil, \prod_{i=k-t+1}^k v_i+\prod_{i=k-t+2}^k v_i \}$, if $v_{k-t}< v_{k-t+1}<2v_{k-t}$ and $t\geq 2$;

\item  $\left\lceil\frac{2 v_{k-1}+2v_k }{3}\right\rceil$, if $v_{k-t}< v_{k-t+1}<2v_{k-t}$ and $t=1$.

\end{enumerate}

\end{theorem}

\begin{table}[!t]
\caption{Lower Bounds on the size of $(\bar{1},2)$-LA} \label{CLBound}
\centering
\begin{tabular}{|c|c|c|}
\hline
Type & Minimum Size & Stimulation Annealing \\
\hline
(2,3,4)  & 16 & 16 \\
(3,3,4)  & 17 & 17 \\
(2,4,4) & 16 & 16 \\
(2,2,3,4)  & 16 & 16 \\
(2,2,5,5)  & 25 & 25 \\
(2,3,3,4)  & 17 & 17 \\
\hline
\end{tabular}
\end{table}

Table \ref{CLBound} presents a lower bound on the size of some certain mixed-level $(\bar{1},2)$-LAs. The first
column lists the types, while the second column displays the lower bound on the size of mixed-level $(\bar{1},2)$-LAs
with the type. The last column presents the size obtained by simulation annealing \cite{T2019}.

A $(\bar{1},t)$-\mbox{LA}$(N; k, (v_1,v_2,\dots, v_k))$ is called {\em optimal} if its size is $(\bar{1},t)$-\mbox{LAN}$(k, (v_1,v_2,\dots, v_k))$.
In what follows, we  will focus on  some constructions for mixed level LAs from combinatorial design theory. Some constructions that produce optimal LAs satisfying the lower bound
in Lemma \ref{case 1} will also be provided.

\section{Constructions of $(\bar{1},t)$-LA$(N;k, (v_1,v_2,\cdots, v_k))$}

Some constructions and existence results for $(\bar{1},t)$-LA$(N;k, (v_1,v_2,\cdots, v_k))$ are presented in this section.

\subsection{A construction for optimal  $(\bar{1},t)$-LA$(\prod_{i=k-t+1}^k v_i;k, (v_1,v_2,\cdots, v_k))$}

Let $2\leq v_1\leq v_2\leq \cdots \leq v_k$. An $N\times k$ array $A$ is called MCA$_2^*(\prod_{i=k-t+1}^k v_i;t,k,(v_1,v_2,\cdots,v_k))$ if
$|\rho(A,T)|=1$ for any $t$-way interaction $T\in {\cal T}=\{\{(k-t+1,v_{k-t+1}),\cdots, (k,v_{k})\}: v_i\in V_i\ (k-t+1\leq i\leq k)\}$ and $|\rho(A,T')|\geq 2$
for any $t$-way interaction $T'\not \in {\cal T}$. If an optimal $(\bar{1},t)$-LA$(N;k, (v_1,v_2,\cdots, v_k))$ with $N=\prod_{i=k-t+1}^k v_i$ exists,
then the following condition must be satisfied.

\begin{lemma} \label{NC}
Let $2\leq v_1\leq v_2\leq \cdots\leq v_{k-t}, 2v_{k-t}\leq v_{k-t+1}\leq v_{k-t+2}\leq \cdots \leq v_k$.
If $A$ is an optimal $(\bar{1},t)$-LA$(N;k,(v_1,v_2,\cdots,v_k))$ with $N=\prod_{i=k-t+1}^k v_i$.
Then, $A$ is an MCA$_2^*(N;t,k,(v_1,v_2,\cdots,v_k))$.
\end{lemma}

\pf Let $A$ be the given optimal $(\bar{1},t)$-LA$(N;k,(v_1,v_2,\cdots,v_k))$ with $N=\prod_{t=k-t+1}^k v_i$. Then, $A$ is an
MCA$(N;t,k,(v_1,v_2,\cdots,v_k)$ by Lemma \ref{(1,t)-LA and MCA}. Because $N=\prod_{t=k-t+1}^k v_i$, we have $|\rho(A,T)|=1$
for any $t$-way interaction $T\in \mathcal{T}$. It follows that
$|\rho(A,T')|\geq 2$ for any $t$-way interaction $T'$ of $A$ from the definition of $(\bar{1},t)$-LA, where $T'\not \in \mathcal{T}$.
Hence, $A$ is an MCA$_2^*(\prod_{i=k-t+1}^k v_i;t,k,(v_1,v_2,\cdots,v_k))$, as desired. \qed

Clearly, an MCA$_2^*(N;t,k,(v_1,v_2,\cdots,v_k))$ is not always a  $(\bar{1},t)$-LA$(N;k,(v_1,v_2,\cdots,v_k))$. Next, we present a
special case of MCA$_2^*$, which produces optimal $(\bar{1},t)$-LAs. First, we introduce the notion of mixed orthogonal arrays (MOAs).

An MOA, or MOA$(N;t,k,(v_1,v_2,\cdots,v_k))$ is an $N\times k$
array with entries in the $i$th column from a set $V_i$ of size $v_i$ such that each $N\times t$
sub-array contains each $t$-tuple occurring an equal number of times as a row. When $v_1=v_2=\cdots=v_k=v$, an MOA is merely an {\em orthogonal array},
denoted by OA$(N;t,k,v)$.

The notion of mixed or asymmetric orthogonal arrays, introduced by Rao \cite{R1973}, have received significant attention in recent years.
These arrays are important in experimental designs as universally optimal fractions of asymmetric factorials. Without loss of generality,
we assume that $v_1\leq v_2\leq \cdots \leq v_k$. By definition of MOA, all $t$-tuples occur in the same number of rows for any $N \times t$
sub-array of an MOA. This number of rows is called  {\em index}. It is obvious that $\binom{k}{t}$ indices exist. We denote it by $\lambda_1,\lambda_2,\cdots, \lambda_{\binom{k}{t}}$. If $\lambda_i\not =\lambda_j$ for
any $i\not =j$, then an MOA is termed as a {\em pairwise distinct index mixed orthogonal array}, denoted by  PDIMOA$(N;t,k,(v_1,v_2,\cdots, v_k))$.
Moreover, if $\lambda_i=1$ for a certain $i\in \{1,2,\cdots, \binom{k}{t}\}$ holds,
then it is termed as PDIMOA$^*(N;t,k,(v_1,v_2,\cdots,v_k))$. It is clear that $N=\prod_{i=k-t+1}^k v_i$ in the definition of PDIMOA$^*$.

\begin{example}
The transpose of the following array is a PDIMOA $^*(24;2,3,(2,4,6))$.
$$
\left(
 \begin{array}{cccccccccccccccccccccccccccccc}
 1 & 0 & 0 & 0 & 1 & 1 &  1 & 0 & 1 & 0 & 0 & 1 & 0 & 1 & 1 & 1 & 0 &
 0 & 0 & 1 & 0 & 1 & 1 & 0\\
 2 & 2 & 2 & 2 & 2 & 2 &  1 & 1 & 1 & 1 & 1 & 1 & 3 & 3 & 3 & 3 & 3 & 3 & 0 & 0 & 0 & 0 & 0 & 0\\
 1 & 2 & 3 & 4 & 5 & 0 &  1 & 2 & 3 & 4 & 5 & 0 & 1 & 2 & 3 & 4 & 5 & 0 & 1 & 2 & 3 & 4 & 5 & 0

\end{array}
\right)
$$  \qed

\end{example}

The following lemma can be easily obtained by the definition of PDIMOA$^*$; therefore, we omit the proof herein.

\begin{lemma}\label{NC1}
Suppose that $v_1 \leq v_2 \leq \cdots  \leq  v_k$. If $A$ is a PDIMOA$^*(\prod_{i=k-t+1}^ k v_i;t,k,(v_1,v_2,\cdots, v_k))$, then
$v_1 <v_2 < \cdots  <  v_k$ and $v_i|v_j$, where $1\leq i\leq k-t$ and $k-t+1\leq j\leq k$.
\end{lemma}

\begin{lemma}\label{PDIMOA-LA}
Let $2< v_1< v_2< \cdots< v_k$. If a PDIMOA$(N;t,k,(v_1,v_2,\dots, v_k))$ exists, then a $(\bar{1},t)$-LA$(N;k,(v_1,v_2,\dots, v_k))$ exists. Moreover, if
$N=\prod_{i=k-t+1}^k v_i$, then the derived $(\bar{1},t)$-LA is  optimal.
\end{lemma}
\pf Let $A$ be a PDIMOA$(N;t,k,(v_1,v_2,\dots, v_k))$. Clearly, $A$ is an MCA. By Lemma \ref{(1,t)-LA and MCA},
we only need to prove that $T_1\not =T_2$ implies $\rho(A,T_1)\not =\rho(A,T_2)$, where $T_1$ and $T_2$ are two $t$-way interactions.
In fact, if $\rho(A,T_1) =\rho(A,T_2)$, then $|\rho(A,T_1)|=|\rho(A,T_2)|$, which contradicts the definition of a PDIMOA.
The optimality can be obtained by Theorem \ref{L-bound}. \qed

We will construct an optimal $(\bar{1},t)$-LA$(N;k,(v_1,v_2,\dots, v_k))$ with $N=\prod_{i=k-t+1}^k v_i$ in terms of PDIMOA$^*$.
First, we have the following simple and useful construction for PDIMOA$^*$. A similar construction for MOAs was
first stated in \cite{CJL2014}.

\begin{construction}\label{f}
Let $b=r_1r_2\cdots r_m<v_2< \cdots<v_k$ and $r_1<r_2<\cdots <r_m$. If a PDIMOA$^*(\prod_{i=k-t+1}^ k v_i;\\t,k,(r_1r_2\cdots r_m,v_2,v_3,\cdots v_k))$ exists,
then a PDIMOA$^*(\prod_{i=k-t+1}^ k v_i;t,k+m-1,(r_1,r_2,\cdots, r_m, v_2,v_3,\cdots, v_k))$ also exists.
\end{construction}

\pf Let $A$ be PDIMOA$^*(N;t,k,(b,v_2,v_3,\cdots v_k))$ with $b=r_1r_2,\cdots r_m$. We can form
an $N\times (k+m-1)$ array $A'$ by replacing the symbols in $V_b$ by those of $V_{r_1}\times V_{r_2} \times \cdots \times V_{r_m}$.
It is easily verified that $A'$ is the required PDIMOA$^*$. \qed

The following construction can be obtained easily; thus, we omit its proof.

\begin{construction}\label{pc2}
Let $a_1<a_2< \cdots<a_k$ and $b_1<b_2< \cdots<b_k$. If both a PDIMOA$^*(\prod_{i=k-t+1}^ k a_i;t,k,\\(a_1,a_2,\cdots, a_k))$ and a PDIMOA$^*(\prod_{i=k-t+1}^ k b_i;t,k, (b_1,b_2,\cdots, b_k))$ exist, then a PDIMOA$^*(\prod_{i=k-t+1}^ k a_ib_i;\\t,k,(a_1b_1,a_2b_2,\cdots, a_kb_k))$ exists. In particular, if both a PDIMOA$^*(\prod_{i=k-t+1}^ k a_i;t,k,(a_1,a_2,\cdots, a_k))$ and an OA$(t,k,v)$ exist, then a PDIMOA$^*(\prod_{i=k-t+1}^ k a_iv^t;t,k,(a_1v,a_2v,\cdots, a_kv))$ exists.
\end{construction}

\subsection{Methods for constructing $(\bar{1},t)$-LA$(N;k, (v_1,v_2,\cdots, v_k))$}

In this subsection, we  modify some constructions for MCAs to the case of $(\bar{1},t)$-LAs.
The next two lemmas provide the ``truncation'' and ``derivation'' constructions, which were first used
to construct mixed CAs.

\begin{lemma}{(Truncation)}\label{trun}
Let $2\leq v_1\leq v_2\leq \cdots \leq v_{i-1}\leq v_i\leq v_{i+1}\leq \cdots \leq v_k$. Then, $(\bar{1},t)$-LAN$(k-1,(v_1,v_2,\cdots,v_{i-1},v_{i+1},\cdots,v_k))\leq$  $(\bar{1},t)$-LAN$(k,(v_1,v_2,\dots,v_{i-1},v_i,v_{i+1},\cdots,v_k))$.
\end{lemma}
\pf Let $A$ be a $(\bar{1},t)$-LA$(N;k,(v_1,v_2,\dots,v_{i-1},v_i,v_{i+1},\cdots,v_k))$ with $N=(\bar{1},t)$-LAN$(k,(v_1,v_2,\dots,\\v_{i-1},v_i,v_{i+1},\cdots,v_k))$.
Delete the $i$th column from $A$ to obtain  a $(\bar{1},t)$-LA$(N;k-1,(v_1,v_2,\cdots,v_{i-1},v_{i+1},\\\cdots,v_k))$.
Thus, $(\bar{1},t)$-LAN$(k-1,(v_1,v_2,\cdots,v_{i-1},v_{i+1},\cdots,v_k))\leq N=(\bar{1},t)$-LAN$(k,(v_1,v_2,\dots,v_{i-1},v_i,\\v_{i+1},\cdots,v_k))$.\qed

\begin{lemma}{(Derivation)}\label{deri}
Let $2\leq v_1\leq v_2\leq \cdots \leq v_{i-1}\leq v_i\leq v_{i+1}\leq \cdots \leq v_k$. Then $v_i \cdot (\bar{1},t-1)$-LAN$(k-1,(v_1,v_2,\cdots,v_{i-1},v_{i+1},\cdots,v_k))\leq$  $(\bar{1},t)$-LAN$(k,(v_1,v_2,\cdots,v_{i-1},v_i,v_{i+1},\cdots,v_k))$, where $t\geq 2$.
\end{lemma}

\pf Let $A$ be a $(\bar{1},t)$-LA$(N;k,(v_1,v_2,\dots, v_k))$ with $N=(\bar{1},t)$-LAN$(k,(v_1,v_2,\dots, v_k))$.
By Lemma \ref{(1,t)-LA and MCA}, $A$ is an MCA and a $(1,t)$-LA. For each $x \in \{0,1,\cdots, v_i-1\}$,
taking the rows in $A$ that involve the symbol $x$ in the $i$th columns and omitting the column yields
an MCA$(N_x; t-1,k-1,(v_1,v_2,\cdots,v_{i-1},v_{i+1},\cdots,v_k))$. We use $A(x)$ to
denote the derived array. Next, we prove that $A(x)$ is a $(1,t-1)$-LA$(N_x;k-1,(v_1,v_2,\cdots,v_{i-1},v_{i+1},\cdots,v_k))$.
In fact, for any $(t-1)$-way interaction $T_1$ and $T_2$ with $T_1\not= T_2$, if $\rho(A(x),T_1)=\rho(A(x),T_2)$, we can form two $t$-way interactions $T_1'$ and $T_2'$  by inserting $(i,x)$ into $T_1$ and $T_2$, respectively.
Hence, $\rho(A,T_1')=\rho(A,T_2')$, where $|\rho(A,T_1')|= |\rho(A(x),T_1)|$ but $T_1'\not=T_2'$. Consequently, $A$ is not a $(1,t)$-LA.
It is clear that $N_i\geq (\bar{1},t-1)$-LAN$(k-1,(v_1,v_2,\cdots,v_{i-1},v_{i+1},\cdots,v_k))$ for $0\leq i\leq v_i-1$. Thus,
$N=N_0+N_1+\cdots+N_{v_i-1}\geq v_i\cdot(\bar{1},t-1)$-LAN$(k-1,(v_1,v_2,\cdots,v_{i-1},v_{i+1},\cdots,v_k))$.  \qed

The following product construction can be used to produce a new LA from old LAs, which is a typical weight
construction in combinatorial design.

\begin{construction}{(Product Construction)}\label{pc1}
If both a $(\bar{1},t)$-LA$(N_1;k,(v_1,v_2,\dots, v_k))$ and an MCA$(N_2;t,k,\\(s_1,s_2,\dots, s_k))$ exist, then
a $(\bar{1},t)$-LA$(N_1N_2;k,(v_1s_1,v_2s_2,\dots, v_ks_k))$ exists. In particular, if both a $(\bar{1},t)$-LA$(N_1;k,(v_1,v_2,\dots, v_k))$ and a $(\bar{1},t)$-LA$(N_2;k,(s_1,s_2,\dots, s_k))$ exist, then a $(\bar{1},t)$-LA$(N_1N_2;k,(v_1s_1,\\v_2s_2,\dots, v_ks_k))$ also exists.
\end{construction}
\pf Let $A=(a_{ij})\ (i \in I_{N_1}, j\in I_k )$ and $B=(b_{ij})\ (i \in I_{N_2}, j\in I_k )$ be the
given $(\bar{1},t)$-LA$(N_1;k,(v_1,v_2,\\\dots, v_k))$  and   MCA$(N_2;t,k,(s_1,s_2,\dots, s_k))$, respectively.
We form an $N_1N_2\times k$ array as follows. For each row $(a_{i1}, a_{i2 }, \cdots, a_{ik})$ of $A$ and each
row $(b_{h1}, b_{h2},\cdots , b_{hk})$ of $B$, include the row
$((a_{i1}, b_{h1}), (a_{i2}, b_{h2}),\\ \cdots, (a_{ik}, b_{hk}))$ as a row of $\overline{A}$, where $1\leq i\leq N_1, 1\leq h\leq N_2$.

From the typical weighting method in design theory, the resultant array $\overline{A}$ is an MCA$(N_1N_2;t, k,(v_1s_1,\\v_2s_2,\dots, v_ks_k))$, as
both $A$ and $B$ are MCAs.  By Lemma \ref{(1,t)-LA and MCA}, we only need to
prove that $\overline{A}$ is a $(1,t)$-LA. Suppose that $\rho(\overline{A},T_1)= \rho(\overline{A},T_2)$, where
$T_1=\{(i,(a_{hi},b_{ci})):i\in I, |I|=t, I\subset \{1,2,\cdots,k\}, h\in I_{N_1},c\in I_{N_2}\}$ and
$T_2=\{(j,(a_{h'j},b_{c'j})):j\in I', |I'|=t, I'\subset \{1,2,\cdots,k\}, h'\in I_{N_1},c'\in I_{N_2}\}$  with $T_1\not =T_2$.
It is noteworthy that the projection on the first component of $T_1$ and $T_2$ is the corresponding $t$-way interaction of $A$,
while the projection on the second component is the corresponding $t$-way interaction of $B$. Therefore,
$A$ is not a $(1,t)$-LA. The first assertion is then proved because a $(\bar{1},t)$-LA$(N_2;k,(s_1,s_2,\dots, s_k))$ is an MCA$(N_2;t,k,(s_1,s_2,\dots, s_k))$.
 The second assertion can be proven by the first assertion.\qed

The following construction can be used to increase the number of levels for a certain factor.

\begin{construction} \label{twice construction}
If a $(\bar{1},t)$-LA$(N;k,(v_1,v_2,\cdots,v_k))$ exists, then a $(\bar{1},t)$-LA$(2N;k,(v_1,v_2,\cdots, v_{i-1}, a, \\ v_{i+1}, \cdots,v_k))$ exists,
where $i\in \{1,2,3,\cdots,k\}$ and $v_i<a\leq 2v_i$.
\end{construction}
\pf Let $A=(a_{ij}), (i\in I_N,j\in I_k)$ be the given  $(\bar{1},t)$-LA$(N;k,(v_1,v_2,\cdots,v_k))$ with entries in the $i$th
column from a set $V_i$ of size $v_i$. For a certain $i\in I_k$, we replace the symbols $0,1,\cdots,a-v_i-1$ in the $i$th column
of $A$ by $v_i,v_{i}+1,\cdots, a-1$, respectively. We denote the resultant array by $A'$. Clearly, permuting the symbols in a
certain column does not affect the property of $(\bar{1},t)$-LAs. Thus, $A'$ is also a $(\bar{1},t)$-LA$(N;k,(v_1,v_2,\cdots,v_k))$, where
 entries in the $i$th column of $A'$ from the set $\{a-v_i, a-v_i+1,\cdots, v_i-1,v_i,v_i+1,\cdots, a-1\}$. Subsequently, write $M=(A^T|(A')^T)^T$.
It is easy to prove that $M$ is a $(1,t)$-LA$(2N;k,(v_1,v_2,\cdots, v_{i-1}, a,v_{i+1}, \cdots,v_k))$ and an MCA$(2N;t,k,(v_1,v_2,\cdots, v_{i-1}, a,v_{i+1}, \cdots,v_k))$.
By Lemma \ref{(1,t)-LA and MCA}, $M$ is the desired array. \qed

The following example illustrates the idea in Construction \ref{twice construction}.

\begin{example}\label {Ex3-1}\rm
The transpose of the following array is a $(\bar{1},2)$-LA$(12;5,(2,2,2,2,3))$
\begin{center}
{\small
 \tabcolsep 1.8pt
 \begin{tabular}{|cccccccccccc|}
 \hline
0 & 0 & 0 & 0 & 0 & 0 & 0 & 1& 1 & 1 & 1 &1 \\
0 & 0 & 0 & 0 & 1 & 1 & 1 & 0& 0 & 1 & 1 &1 \\
0 & 0 & 0 & 1 & 0 & 1 & 1 & 1& 1 & 0 & 0 &1 \\
0 & 0 & 1 & 0 & 1 & 0 & 1 & 1& 1 & 0 & 1 &0 \\
0 & 2 & 1 & 0 & 1 & 1 & 2 & 0& 1 & 0 & 2 &2 \\
\hline
\end{tabular}
}
\end{center}

\noindent
Replace the symbols $0,1$ by $2,3$ in the $3$th column, respectively. Juxtapose  two such arrays from top to bottom to obtain the following array $M$; we list it
as its transpose to conserve space.
\begin{center}
{\small
 \tabcolsep 1.8pt
 \begin{tabular}{|cccccccccccc|cccccccccccc|}
 \hline
0 & 0 & 0 & 0 & 0 & 0 & 0 & 1& 1 & 1 & 1 &1 &0 & 0 & 0 & 0 & 0 & 0 & 0 & 1& 1 & 1 & 1 &1  \\
0 & 0 & 0 & 0 & 1 & 1 & 1 & 0& 0 & 1 & 1 &1 &0 & 0 & 0 & 0 & 1 & 1 & 1 & 0& 0 & 1 & 1 &1  \\
0 & 0 & 0 & 1 & 0 & 1 & 1 & 1& 1 & 0 & 0 &1 &2 & 2 & 2 & 3 & 2 & 3 & 3 & 3& 3 & 2 & 2 &3  \\
0 & 0 & 1 & 0 & 1 & 0 & 1 & 1& 1 & 0 & 1 &0 &0 & 0 & 1 & 0 & 1 & 0 & 1 & 1& 1 & 0 & 1 &0  \\
0 & 2 & 1 & 0 & 1 & 1 & 2 & 0& 1 & 0 & 2 &2 &0 & 2 & 1 & 0 & 1 & 1 & 2 & 0& 1 & 0 & 2 &2  \\
\hline
\end{tabular}
}
\end{center}
\noindent
It is easy to verify that $M$ is a $(\bar{1},2)$-LA$(24;5,(2,2,4,2,3))$.

Replace the symbol $0$ by $2$ in the $3$th column. Juxtapose  two such arrays from top to bottom to obtain the following array $M'$; we list it
as its transpose to conserve space.
\begin{center}
{\small
 \tabcolsep 1.8pt
 \begin{tabular}{|cccccccccccc|cccccccccccc|}
 \hline
0 & 0 & 0 & 0 & 0 & 0 & 0 & 1& 1 & 1 & 1 &1 &0 & 0 & 0 & 0 & 0 & 0 & 0 & 1& 1 & 1 & 1 &1  \\
0 & 0 & 0 & 0 & 1 & 1 & 1 & 0& 0 & 1 & 1 &1 &0 & 0 & 0 & 0 & 1 & 1 & 1 & 0& 0 & 1 & 1 &1  \\
0 & 0 & 0 & 1 & 0 & 1 & 1 & 1& 1 & 0 & 0 &1 &2 & 2 & 2 & 1 & 2 & 1 & 1 & 1& 1 & 2 & 2 &1  \\
0 & 0 & 1 & 0 & 1 & 0 & 1 & 1& 1 & 0 & 1 &0 &0 & 0 & 1 & 0 & 1 & 0 & 1 & 1& 1 & 0 & 1 &0  \\
0 & 2 & 1 & 0 & 1 & 1 & 2 & 0& 1 & 0 & 2 &2 &0 & 2 & 1 & 0 & 1 & 1 & 2 & 0& 1 & 0 & 2 &2  \\
\hline
\end{tabular}
}
\end{center}

\noindent
It is easy to verify that $M'$ is a $(\bar{1},2)$-LA$(24;5,(2,2,3,2,3))$.\qed
\end{example}

 {\bf Remark:} Construction \ref{twice construction} may produce an optimal $(\bar{1},t)$-LA. For example, a $(\bar{1},2)$-LA$(16;(2,2,3,4))$ is shown in
Table 1. By Construction \ref{twice construction}, we can obtain a $(\bar{1},2)$-LA$(32;(2,2,3,8))$, which is optimal by Lemma \ref{case 3}.

Fusion is an effective construction for MCAs from CAs. It causes any $d\ge 2$ levels to be identical;
for example, see \cite{CGRS2010}. As with CAs, fusion for $(\bar{1},t)$-LAs guarantees
the extension of uniform constructions to mixed cases. However, fusion for a $(\bar{1},t)$-LA$(N;k,v)$ may not
produce mixed-level $(\bar{1},t)$-LAs. This problem can be circumvented by introducing the notion of detecting arrays (DAs). If, for
any ${\cal T}\subseteq \mathcal {I}_t$ with $|{\cal T}| = d$ and any $T\in \mathcal {I}_t$, we
have $\rho(A, T)\subseteq \rho(A, {\cal T}) \Leftrightarrow T\in {\cal T},$ then the array $A$
is called a $(d,t)$-DA or a $(d,t)$-DA$(N; k,v)$.

\begin{construction}{(Fusion)}\label{fusion}
Suppose that $A$ is a $(1,t)$-DA$(N; k,v)$ with $t\geq 2$. If $A$ is also
a $(\lceil \frac{v}{v_i}\rceil,t)$-LA$(N; k,v)$, then a $(\bar{1},t)$-LA$(N;k,(v,\cdots,v,v_i,v,\cdots, v))$
exists, where $ 2\leq v_i< v$.
\end{construction}

\pf Let $A$ be a $(1,t)$-DA$(N; k,v)$ over the symbol set $V$ of size $v$.
Let $a_1+a_2+\cdots+a_{v_i}=v$, where $a_i(i=1,2,\cdots,v_i)\geq 1$. We can select one $a_i$
such that $a_i=\lceil \frac{v}{v_i}\rceil$ and $a_i\geq a_j$, where $1\leq i\not =j\leq v_i$.
We select $a_1, a_2,\cdots, a_{v_i}$  elements from $V$ in the $i$th column of $A$ to form the element sets $A_i(1\leq i\leq v_i)$, respectively.
The elements in $A_i(1\leq i\leq v_i)$ are  identical with $1,2,\cdots, v_i$, respectively.
Then, we obtain an $N\times k$ array $A'$. Clearly, $A'$ is an MCA. We only need to prove that $A'$ is a $(1,t)$-LA by Lemma \ref{(1,t)-LA and MCA}, i.e.,
for any two distinct $t$-way interactions $T_1=\{(a_1,u_{a_1}),\cdots,(a_t,u_{a_t})\}$ and $T_2=\{(b_1,s_{b_1}),  \cdots,(b_t,s_{b_t})\}$,
we have $\rho(A',T_1)\not =\rho(A',T_2)$. It is clear that $\rho(A,T_1)=\rho(A',T_1)$ and $\rho(A',T_2)=\rho(A,T_2)$ when $i\not \in \{a_1,\cdots,a_t\}$ and $i\not \in \{b_1,\cdots,b_t\}$. Hence, $\rho(A',T_1)\not=\rho(A',T_2)$.

When $i\in  \{a_1,\cdots,a_t\}$ and $i\not \in \{b_1,\cdots,b_t\}$, we can obtain a $t$-way interaction $T_1'=\{(a_1,u_{a_1},\cdots,(i,a),\\\cdots,(a_t,u_{a_t})\}$ of $A$,
where $a\in A_{u_i}$. If $\rho(A',T_1) =\rho(A',T_2)$, then
$\rho(A,T_1')\subset \rho(A',T_1)=\rho(A',T_2)=\rho(A,T_2)$. However, $T_1'\not =T_2$; as such, it is a contradiction that $A$ is a $(1,t)$-DA$(N;k,v)$.
If $i\not \in  \{a_1,\cdots,a_t\}$ and $i \in \{b_1,\cdots,b_t\}$, then the similar argument can prove the conclusion.

When $i\in  \{a_1,\cdots,a_t\}$ and $i\in \{b_1,\cdots,b_t\}$, it is clear that $\rho(A',T_1)\not =\rho(A',T_2)$ if $u_i\not =s_i$.
The case $u_i=s_i$ remains to be considered. Without loss of generality, suppose that $a_j$ elements are identical with $u_i$.
It is clear that $T_1$ and $T_2$ can be obtained from ${\cal T}_1$ and ${\cal T}_2$ by fusion, respectively, where ${\cal T}_1$ and ${\cal T}_2$
are sets of $t$-way interactions with $|{\cal T}_1|=|{\cal T}_2|=a_j$. If $\rho(A',T_1) =\rho(A',T_2)$,
then $\rho(A',T_1)=\rho(A,{\cal T}_1) =\rho(A',T_2)=\rho(A,{\cal T}_2)$. It is a contradiction that $A$
is a $(\lceil \frac{v}{v_i}\rceil,t)$-LA$(N;k,v)$ because the existence of $(\lceil \frac{v}{v_i}\rceil,t)$-LA$(N;k,v)$ implies the existence of
$(a_j,t)$-LA$(N;k,v)$ \cite{CM2008}.\qed

Constructions \ref{twice construction} and \ref{fusion} provide an effective and efficient method
to construct a mixed-level $(\bar{1},t)$-LA from a $(1,t)$-LA$(N;k,v)$.
The existence of $(d,t)$-DA$(N;k,v)$ with $d\geq 1$ implies the existence of $(d,t)$-LA$(N;k,v)$ \cite{CM2008}.
Hence, the array $A$ in Construction \ref{fusion} can be obtained by a $(d,t)$-DA$(N;k,v)$, which
is characterized in terms of super-simple OAs. The existence of super-simple OAs can be found
in \cite{Chen2011,H2000,STY2012,SW2016,SY2014, TY2011}. It is noteworthy that the derived array
is not optimal. In the remainder of this section, we present two ``Roux-type'' recursive constructions\cite{R1987}.

\begin{construction} \label{increasing one group}
If both a $(\bar{1},t)$-LA$(N_1;k,(v_1,v_2,\cdots,v_k))$ and a $(\bar{1},t-1)$-LA$(N_2;k-1,(v_1,v_2,\cdots,\\v_{i-1},v_{i+1},\cdots,v_k))$ exist,
then a $(\bar{1},t)$-LA$(N_1+eN_2;k,(v_1,v_2,\cdots, v_{i-1}, v_i+e, v_{i+1}, v_{i+2}\cdots,v_k))$ exists, where $e\geq 0$.
\end{construction}

\pf Let $A$ and $B$ be the given $(\bar{1},t)$-LA$(N_1;k,(v_1,v_2,\cdots,v_k))$ and $(\bar{1},t-1)$-LA$(N_2;k-1,(v_1,v_2,\cdots,\\v_{i-1},v_{i+1},\cdots,v_k))$, respectively. Clearly, if $e=0$, then $A$ is the required array. Now, suppose that $e\geq 1$. Insert a column vector $(j,j,\cdots,j)$ of length $N_2$ to the front of the $i$th
column of $B$ to form an $N_2\times k$ array $B_j$, where $j\in \{v_i,v_i+1,v_i+2,\cdots,v_i+e-1\}$. Let $M=(A^T|B_{v_i}^T|B_{v_i+1}^T|\cdots|B_{v_i+e-1}^T)^T$. Clearly,
$M$ is an MCA$(N_1+eN_2;t,k,(v_1,v_2,\cdots, v_{i-1}, v_i+e, v_{i+1}, v_{i+2}\cdots,v_k))$ \cite{CSWY2011}. By Lemma \ref{(1,t)-LA and MCA}, we only need to
prove that $M$ is a $(1,t)$-LA, i.e., $\rho(M,T_1)\not =\rho(M,T_2)$ for any two distinct $t$-way interactions $T_1$ and $T_2$, where
$T_1=\{(a_1,u_{a_1}),\cdots,(a_t,u_{a_t})\}$ and $T_2=\{(b_1,s_{b_1}),  \cdots,(b_t,s_{b_t})\}$. Next, we distinguish the following cases.

\noindent{Case 1. }  $i\not \in \{a_1,\cdots, a_t\}$ and $i\not \in \{b_1,\cdots, b_t\}$

In this case, because $A$ is a $(\bar{1},t)$-LA, $\rho(A,T_1)\not= \rho(A,T_2)$, $\rho(M,T_1)\not= \rho(M,T_2)$ as $A$ is part of $M$.

\noindent{Case 2. }  $i\not \in \{a_1,\cdots, a_t\}$ and $i\in \{b_1,\cdots, b_t\}$ or $i \in \{a_1,\cdots, a_t\}$ and $i \not \in \{b_1,\cdots, b_t\}$

When $i\not \in \{a_1,\cdots, a_t\}$ and $i\in \{b_1,\cdots, b_t\}$, if $s_i\not \in \{v_i,v_i+1,\cdots, v_i+e-1\}$,
then $\rho(A,T_1)\not= \rho(A,T_2)$. Thus, $\rho(M,T_1)\not= \rho(M,T_2)$. If $s_i \in \{v_i,v_i+1,\cdots, v_i+e-1\}$,
then $T_2$ must be included by rows of $B_i$, where $i\in \{v_i,v_i+1,\cdots, v_i+e-1\}$; however, it must not be included by any row of $A$.
Clearly, $T_1$ must be included by some rows of $A$. Consequently, $\rho(M,T_1)\not= \rho(M,T_2)$.
When $i \in \{a_1,\cdots, a_t\}$ and $i \not \in \{b_1,\cdots, b_t\}$, the same argument can prove the conclusion.

\noindent{Case 3. }   $i \in \{a_1,\cdots, a_t\}$ and $i\in \{b_1,\cdots, b_t\}$

Clearly, $\rho(M,T_1)\not= \rho(M,T_2)$ holds whenever $u_i\not =s_i$.
If $u_i=s_i\not \in \{v_i,v_i+1,\cdots,v_i+e-1\}$, then $\rho(A,T_1)\not= \rho(A,T_2)$, which
implies that $\rho(M,T_1)\not= \rho(M,T_2)$. If $u_i=s_i\in \{v_i,v_i+1,\cdots,v_i+e-1\}$, then $T_1$ and $T_2$ must
be included by some rows for a certain $B_i$, where $i\in \{v_i,v_i+1,\cdots,v_i+e-1\}$. Because $B$ is a $(\bar{1},t-1)$-LA,
$\rho(B_i, T_1)\not =\rho(B_i,T_2)$, which implies $\rho(M,T_1)\not= \rho(M,T_2)$. \qed

More generally, we have the following construction.

\begin{construction}\label {increasing two groups}
Let $p\geq 0, q\geq 0$ and $1\leq i< j\leq k$.
If a $(\bar{1},t)$-LA$(N_1;k,(v_1,v_2,\cdots,v_{i-1},v_{i},v_{i+1},\cdots,\\ v_{j-1},v_{j},v_{j+1},\cdots,v_k))$,
$(\bar{1},t-1)$-LA$(N_2;k-1,(v_1,v_2,\cdots ,v_{i-1},v_{i+1},\cdots,v_k))$, a $(\bar{1},t-1)$-LA$(N_3;k-1,(v_1,v_2,\cdots ,v_{j-1},v_{j+1},\cdots,v_k))$
and $(\bar{1},t-2)$-LA$(N_4;k-2,(v_1,v_2,\cdots,v_{i-1},v_{i+1},\cdots,v_{j-1},v_{j+1},\cdots,\\v_k))$ exist,
then a $(\bar{1},t)$-LA$(N;k,(v_1,v_2,\cdots,v_{i-1},v_{i}+p,v_{i+1},\cdots, v_{j-1},v_{j}+q,v_{j+1},\cdots,v_k))$ exists, where  $N=N_1+pN_2+qN_3+pqN_4$.

\end{construction}
\pf
We begin with a $(\bar{1},t)$-LA$(N_1;k,(v_1,v_2,\cdots,v_{i-1},v_{i},v_{i+1},\cdots, v_{j-1},v_{j},v_{j+1},\cdots,v_k))$,
an $N_1\times k$ array $A$ that is on $V_1\times\cdots  \times V_{i-1}\times V_{i}'\times V_{i+1}\times\cdots \times V_{j-1}\times V_{j}'\times V_{j+1}\times\cdots \times
V_k$. Let $H_1$ and $H_2$ be two sets with $|H_1|=p$ and $|H_2|=q$ such that $H_1\bigcap V_i'=\emptyset$ and $H_2\bigcap V_j'=\emptyset$, respectively.
Suppose that $B'$, an $N_2\times (k-1)$ array, is a $(\bar{1},t-1)$-LA$(N_2;k-1,(v_1,v_2,\cdots ,v_{i-1},v_{i+1},\cdots,v_k))$,
which is on $V_1\times\cdots  \times V_{i-1}\times V_{i+1}\times\cdots \times V_k$. For each row $(a_1$, $a_2,\cdots, a_{i-1},a_{i+1},\cdots, a_k)$ of $B'$, add $x\in H_1$ to obtain a $k$-tuple $(a_1$, $a_2,\cdots, a_{i-1},x,a_{i+1},\cdots,a_k)$. Then, we obtain a $pN_2\times k$ array from $B'$, denoted by $B$. Similarly,
from a $(\bar{1},t-1)$-LA$(N_3;k-1,(v_1,v_2,\cdots ,v_{j-1},v_{j+1},\cdots,v_k))$, we obtain a $qN_3\times k$ array, denoted by $C$.
For each pair $(x,y)\in H_1\times H_2$, we construct $k$-tuple $(a_1, a_2,\cdots, a_{i-1},x,a_{i+1},\cdots,a_{j-1},y,a_{j+1},\cdots,a_k)$ for each row of the given
$(\bar{1},t-2)$-LA$(N_4;k-2,(v_1,v_2,\cdots,v_{i-1},v_{i+1},\cdots,v_{j-1},v_{j+1},\cdots,v_k))$. These tuples result in a $pqN_4\times k$ array, denoted by $D$.

 Denote $V_i'\cup H_1=V_i$, $V_j'\cup H_2=V_j$ and
$
 F=\left(
 \begin{array}{l}
 A\\
 B\\
 C\\
 D\\
 \end{array}
 \right)
 $.  We claim that $F$, an $(N_1+pN_2+qN_3+pqN_4)\times k$ array, is a $(\bar{1},t)$-LA$(N;k,(v_1,v_2,\cdots,v_{i-1},v_{i}+p,v_{i+1},\cdots, v_{j-1},v_{j}+q,v_{j+1},\cdots,v_k)$
 which is on $V_1\times\cdots  \times V_{i-1}\times V_{i}\times V_{i+1}\times\cdots \times V_{j-1}\times V_{j}\times V_{j+1}\times\cdots \times
 V_k$.

Clearly, $F$ is an MCA$(N;t,k,(v_1,v_2,\cdots,v_{i-1},v_{i}+p,v_{i+1},\cdots, v_{j-1},v_{j}+q,v_{j+1},\cdots,v_k)$.
To prove this assertion, we only need to demonstrate that $\rho(F,T_a)\not =\rho(F,T_b)$ for any two distinct $t$-way interactions $T_a=\{(a_1,u_{a_1}),\cdots,(a_t,u_{a_t})\}$ and $T_b=\{(b_1,v_{b_1}),\cdots,(b_t,v_{b_t})\}$. By similar argument as the proof of Construction \ref{increasing one group}, we can prove the conclusion except for the case
where $i,j\in \{a_1,a_2,\cdots,a_t\}$ and $i,j\in \{b_1,b_2,\cdots,b_t\}$, $u_i=v_i\in H_1$, and $ u_j=v_j\in H_2$. In this case, $T_a$ and $T_b$ are only included by some
rows of $D$. If $\rho(F,T_a)=\rho(F,T_b)$, then $\rho(D,T_a) =\rho(D,T_b)=\rho(F,T_a)=\rho(F,T_b)$. Consequently, $\rho(D,T_a\setminus \{(i,u_i),(j,u_j)\}) =\rho(D,T_b\setminus \{(i,u_i),(j,u_j)\})$, which implies that  $\rho(D',T_a\setminus \{(i,u_i),(j,u_j)\}) =\rho(D',T_b\setminus \{(i,u_i),(j,u_j)\})$
by the construction of $D$. It is a contradiction with $D'$ being a $(\bar{1},t-2)$-LA$(N_4;k-2,(v_1,v_2,\cdots,v_{i-1},v_{i+1},\cdots,v_{j-1},v_{j+1},\cdots,v_k))$. The proof is completed.\qed

\subsection{Optimal $(\bar{1},t)$-LA$(\prod_{i=k-t+1}^k;k, (v_1,v_2,\cdots, v_k))$}

In this subsection, some series of optimal mixed-level $(\bar{1},t)$-LAs are presented. First, we list some known results for later use.

\begin{lemma}\label{OA(t,t+1,v)} \cite{HSS1999}
An OA$(v^t;t,t+1,v)$ exists for any integer $v\geq 2, t\geq  2$.
\end{lemma}

The existence of PDIMOA$^*(t,t+1,(v_1,v_2,\cdots, v_t))'$s is determined completely by the following theorem.

\begin{theorem}\label{PDIMOA(t,t+1)}
Let $v_1<v_2<\cdots <v_{t+1}$. A PDIMOA$^*(\prod_{i=2}^{t+1} v_i;t,t+1,(v_1,v_2,\cdots, v_t,v_{t+1}))$ exists if and only if $v_1|v_i$ for
$2\leq i\leq t+1$.
\end{theorem}
\pf The necessity can be easily obtained by Lemma \ref{NC1}. For sufficiency, we write $v_i=v_1r_i$ for $i=2,3,\cdots, t+1$.
Clearly, $r_i\geq 2$ and $r_i\not= r_j$ for $2\leq i\not=j\leq t+1$.  We list all $t$-tuples from $Z_{r_2}\times Z_{r_3}\times \cdots \times Z_{r_{t+1}}$ to form an MOA$(\prod_{i=2}^{t+1}
r_i;t,t,(r_2,r_3,\cdots, r_t, r_{t+1})$, which is also a PDIMOA$^*(\prod_{i=2}^{t+1}r_i;t,t+1,(1,r_2,r_3,\cdots, r_t, r_{t+1})$. Apply Construction \ref{pc2} with
an OA$(v_1^t;t,t+1,v_1)$ given by Lemma \ref{OA(t,t+1,v)} to obtain the required PDIMOA$^*$. \qed

More generally, we have the following results.
\begin{theorem}
Let $v_1<v_2<\cdots<v_k$ and $v_i=k_iv_1v_2\cdots v_{k-t}$, where $k_i\geq 2$, $i=k-t+1,k-t+2,\cdots, k$. Then, a PDIMOA$^*(\prod _{i=k-t+1}^ k v_i; t,k,(v_1,v_2,\cdots, v_k))$ exists.
\end{theorem}


\pf Let $M=v_1v_2\cdots v_{k-t}$. Then, $v_i=Mk_i$, where $i=k-t+1,\cdots, k$.  By Theorem \ref{PDIMOA(t,t+1)}, a PDIMOA$^*(N;t,t+1,(M,v_{k-t+1},\cdots,v_k))$ with $N=\prod _{i=k-t+1}^ k v_i$ exists. Apply Construction \ref{f} to obtain a PDIMOA$^*(\prod _{i=k-t+1}^ k v_i, t,k,(v_1,v_2,\cdots, v_k))$ as desired. \qed

\begin{theorem}\label{LA(2,3,v)}
Let $v_1\leq v_2\leq v_3$ with $v_2\geq 2v_1$. Then, an optimal $(\bar{1},2)$-LA$(v_2v_3;3,(v_1,v_2,v_3))$ exists.
\end{theorem}

\pf First, we construct a $v_2v_3\times 3$ array $A=(a_{ij})$ : $a_{i+rv_3,1}=(i-1+r)\%v_1$, where $i=1,2,\cdots, v_3$ and $r=0,1,\cdots,v_2-1$;
$a_{i,2}=\left \lfloor\frac{i-1} {v_3} \right\rfloor $ and $a_{i,3}=(i-1)\%v_3$ for $i=1,2,\cdots, v_2v_3$.

We will prove that $A$ is an optimal $(\bar{1},2)$-LA. Optimality is guaranteed by Theorem \ref{L-bound}.
It is clear that $A$ is MCA$^*_2(v_2v_3, (v_1,v_2,v_3))$.  Consequently, $|\rho(A,\{(1,a),(2,b)\})|\geq 2,
|\rho(A,\{(1,c),(3,d)\})|\geq 2$ and $|\rho(A,\{(2,e),(3,f)\})|=1$, where $a,c\in V_1,b,e\in V_2, d,f\in V_3$.
It is clear that $\rho(A,\{(1,a),(2,b)\})\not =\rho(A,\{(2,e),(3,f)\})$ and $\rho(A,\{(1,c),(3,d)\})\not =\rho(A,\{(2,e),(3,f)\})$.
We only need to prove $\rho(A,\{(1,a),(2,b)\}) \\\not=\rho(A,\{(1,c),(3,d)\})$. In fact, by construction,
$\rho(A,\{(1,a),(2,b)\}) \subset \{rv_3+1,rv_3+2,\cdots, (r+1)v_3\}$
for a certain $r\in \{0,1,2,\cdots, v_2-1\}$ but $\{i,i+v_1v_3\}\subset \rho(A,\{(1,c),(3,d)\})$, where $i\in \{1,2,\cdots, v_1v_3\}$,
which implies $\rho(A,\{(1,a),(2,b)\})\not =\rho(A,\{(1,c),(3,d)\})$.  Thus, $A$ is a $(\bar{1},t)$-LA by Lemma \ref{(1,t)-LA and MCA}.    \qed

The following example illustrates the idea in Theorem \ref{LA(2,3,v)}.

\begin{example}
The transpose of the following array is an optimal $(\bar{1},2)$-LA$(42;3,(3,6,7))$
\begin{center}
{\small
 \tabcolsep 2.0pt
\begin{tabular}{|ccccccc|ccccccc|ccccccc|ccccccc|ccccccc|ccccccc|}
\hline
0 & 1 &2 & 0& 1 &2 &0 & 1 & 2  & 0& 1 &2 & 0 &1 & 2 & 0 &1 & 2& 0 &1&2 & 0 & 1 &2 & 0& 1 &2 &0 & 1 & 2  & 0& 1 &2 & 0 &1 & 2 & 0 &1 & 2& 0 &1&2 \\
0 & 0 &0 & 0& 0 &0 &0 & 1 & 1  & 1& 1 &1 & 1 &1 & 2 & 2 &2 & 2& 2 &2&2 & 3 & 3 &3 & 3& 3 &3 & 3 & 4 & 4 &4 & 4& 4 &4& 4 & 5 & 5 &5 & 5& 5 &5& 5 \\
0 & 1 &2 & 3& 4 &5 &6 & 0 & 1 &2 & 3& 4 &5 &6 & 0 & 1 &2 & 3& 4 &5 &6& 0 & 1 &2 & 3& 4 &5 &6& 0 & 1 &2 & 3& 4 &5 &6& 0 & 1 &2 & 3& 4 &5 &6 \\
\hline
\end{tabular}
}
\end{center}
\qed
\end{example}

\begin{theorem}\label{LA(1,K,v)}
Let $2\leq w<v$ with $v\geq 2w$. Then, an optimal $(\bar{1},1)$-LA$(v;w+1,(w,w,\cdots,w,v))$ exists.
\end{theorem}
\pf
First, we construct a $2w\times (w+1)$ array $A=(a_{ij})$ as follows:
$$
A=\left[
\begin{array}{ccccc}

0 & 0  & \cdots & 0 &  0 \\
1 & 1  & \cdots & 1 &  1 \\
\vdots & \vdots  &\vdots  & \vdots &\vdots \\
w-1 & w-1  &\cdots & w-1 &w-1 \\
0 & 1  & \cdots & w-1 &  w\\
1 & 2  & \cdots & 0 &   w+1\\
\vdots & \vdots  &\vdots  & \vdots   \\
w-1 & 0  &\cdots  & w-2 & 2w-1 \\

\end{array}
\right]
$$

When $v>2w$, let $C=(c_{ij})$ be a $(v-2w)\times (w+1)$ array with $c_{i,(w+1)}=i-1$ for $i=2w+1,2w+2,\cdots, v$ and $c_{i,j}$ be an arbitrary element for
$\{0,1,\cdots, w-1\}$ with $i=2w+1,2w+2,\cdots,v, j=1,2,\cdots,w$. Let $M=A$ and $N=(A^T|C^T)^T$. It is easy to prove that $M$ and $N$ are the required arrays
if $v=2w$ and $v>2w$, respectively. \qed

\section{Concluding Remarks}
LAs can be used to generate test suites for combinatorial testing and identify interaction faults in component-based systems.
In this study, a lower bound on the size of $(\bar{1},t)$-LAs with mixed levels was determined.
In addition, some constructions of $(\bar{1},t)$-LAs were proposed. Some of these constructions produce optimal locating arrays.
Based on the constructions, some infinite series of optimal locating arrays satisfying the lower bound in Lemma \ref{case 1} were presented.
Obtaining new constructions for mixed-level $(\bar{1},t)$-LAs and providing more existence results are potential future directions.

\end{document}